\documentclass[english]{paper}

\title{On Quantifiers for Quantitative Reasoning}
\author{Matteo Capucci\\\vspace*{-2ex}{\normalsize MSP Group, University of Strathclyde, Glasgow (UK)}\\\small\texttt{\href{mailto:matteo.capucci@gmail.com}{matteo.capucci@gmail.com}}}
\date{\today}

\usepackage{fonttable}
\usepackage{wasysym}
\usepackage{expl3}
\usepackage{xparse, xpatch}

\usepackage{todonotes}
\usepackage{tabularx, multirow}
\usepackage{stackrel, stackengine, old-arrows}

\usepackage{mathrsfs, amssymb, amsfonts, stmaryrd}

\usepackage{kpfonts}
\usepackage{newpxtext}
\usepackage{inconsolata}

\usepackage{pgfplots}

\usepackage{quiver}

\newenvironment{eqalign}{\begin{equation}\begin{aligned}}{\end{aligned}\end{equation}}
\newenvironment{eqalign*}{\begin{equation*}\begin{aligned}}{\end{aligned}\end{equation*}}

\tikzset{
  relation/.style={
    draw=none,
    every to/.append style={
      edge node={node [sloped, allow upside down, auto=false]{$#1$}}}
  }
}

\tikzcdset{
  mono/.code={
    \pgfsetarrows{tikzcd to reversed-tikzcd to}
  }
}
\tikzcdset{
  epi/.code={
    \pgfsetarrowsend{tikzcd double to}
  }
}
\tikzcdset{
  into/.code={
    \pgfsetarrows{tikzcd right hook-tikzcd to}
  }
}
\tikzcdset{
  twocell/.style={Rightarrow, shorten >= 3ex, shorten <= 3ex}
}

\tikzcdset{
  row sep/normal={
    6ex
  },
  column sep/normal={
    8ex
  }
}

\newcommand{\Overset}[2]{%
  \mathop{#2}\limits^{\vbox to -.1ex{%
  \kern -1.8ex\hbox{$#1$}\vss}}%
}
\newcommand{\Underset}[2]{%
  \mathop{#2}\limits_{\vbox to .1ex{%
  \kern -.6ex\hbox{$#1$}\vss}}%
}

\mathchardef\dash="2D

\newcommand{\e}{\mathrm{e}}

\renewcommand{\exp}{\operatorname{exp}}

\DeclareMathOperator{\argmax}{\mathrm{argmax}}

\newcommand{\longto}{\longrightarrow}

\newcommand{\nisoto}[1]{\xrightarrow[\sim]{#1}}
\newcommand{\nisofrom}[1]{\xleftarrow[\sim]{#1}}
\newcommand{\isofrom}[1]{\nisofrom{}}
\newcommand{\isoto}[1]{\nisoto{}}

\newcommand{\const}{\mathsf{const}}

\newcommand{\de}{\mathrm{d}}

\newcommand{\conv}[1][]{\underset{{#1}}{\longrightarrow}}

\newcommand{\biglens}[2]{
	 \begin{pmatrix}{\vphantom{f_f^f}#1} \\ {\vphantom{f_f^f}#2} \end{pmatrix}
}
\newcommand{\littlelens}[2]{
	 \begin{psmallmatrix}{\vphantom{f}#1} \\ {\vphantom{f}#2} \end{psmallmatrix}
}
\newcommand{\lens}[2]{
  \relax\if@display
	 \biglens{#1}{#2}
  \else
	 \littlelens{#1}{#2}
  \fi
}

\usepackage{xstring}
\newcommand{\cat}[1]{
  \relax
  \StrLen{#1}[\catarglen]
  \ifnum\catarglen=1
    \mathcal{#1}
  \else
    \mathbf{#1}
  \fi
}
\newcommand{\dblcat}[1]{\cat{\mathbb #1}}

\newcommand{\true}{\mathsf{true}}

\newcommand{\iso}[1][]{\overset{#1}{\cong}}

\newcommand{\adj}{\dashv}

\newcommand{\Cat}{\cat{Cat}}

\newcommand{\Meas}{\cat{Meas}}

\newcommand{\op}{\mathsf{op}}

\newcommand{\VCat}[1]{{#1}\dash\Cat}

\newcommand{\psum}[2][p]{\mathchoice{{\bigoplus_{#2}}^{#1}}{\bigoplus_{#2}^{#1}}{\bigoplus_{#2}^{#1}}{\bigoplus_{#2}^{#1}}}
\newcommand{\hsum}[1]{\psum[*]{#1}}

\newcommand{\pmean}[2][p]{\int^{#1}_{#2}}
\newcommand{\phmean}[2][-p]{\pmean[#1]{#2}}
\newcommand{\hmean}[1]{\pmean[*]{#1}}

\newcommand{\pexists}[2][p]{{\textstyle\pmean[{#1}]{#2}}}
\newcommand{\pforall}[2][-p]{{\textstyle\phmean[{#1}]{#2}}}

\newcommand{\PosMulReals}{[1,\infty]_\tensor}
\newcommand{\NegMulReals}{[0,1]_\tensor}
\newcommand{\PosAddReals}{[0,\infty]_\add}
\newcommand{\NegAddReals}{[-\infty,0]_\add}
\newcommand{\PosReals}{[0,\infty]}
\newcommand{\MulReals}{\PosReals_\tensor}
\newcommand{\AddReals}{[-\infty,\infty]_\add}

\newcommand{\tensor}{\otimes}
\newcommand{\cotensor}{\mathbin{\tensor^*}}
\newcommand{\add}{\oplus}
\newcommand{\coadd}{\mathbin{\add^*}}
\newcommand{\with}{\coadd}
\newcommand{\Zero}{\mathbf{0}}
\newcommand{\One}{\mathbf{1}}
\renewcommand{\true}{{\bf true}}
\newcommand{\false}{{\bf false}}

\newcommand{\Rel}{\cat{Rel}}

\DeclareMathOperator{\softmax}{softmax}

\newcommand{\psoftmax}[1][p]{{#1}\dash\softmax}
\renewcommand{\argmax}{\operatorname{argmax}}

\newcommand{\LT}{{\bf LT}}

\newcommand{\PL}{\mathbb{PL}}
\newcommand{\AL}{\mathbb{AL}}
\newcommand{\QPL}{\mathbb{QPL}}

\newcommand{\M}{{\cal M}}
\newcommand{\Quant}{\dblcat{Qnt}}

\newcommand{\analogous}{\ \begin{tikzcd}[ampersand replacement=\&]{}\arrow[squiggly,tail reversed,r]\&{}\end{tikzcd}\ }

\newcommand{\entails}{\vdash}
\newcommand{\Prop}{{\bf Prop}}

\newcommand{\sem}[1]{\llbracket{#1}\rrbracket}

\newcommand{\esssup}{\operatorname{ess\,sup}}
\newcommand{\essinf}{\operatorname{ess\,inf}}

\allowdisplaybreaks
\begin{document}
	\maketitle
	\begin{abstract}
		We explore a kind of first-order predicate logic with intended semantics in the reals.
		Compared to other approaches in the literature, we work predominantly in the multiplicative reals $[0,\infty]$, showing they support three generations of connectives, that we call \emph{nonlinear}, \emph{linear additive} and \emph{linear multiplicative}.
		Means and harmonic means emerge as natural candidates for bounded existential and universal quantifiers, and in fact we see they behave as expected in relation to the other logical connectives.
		We explain this fact through the well-known fact that min/max and arithmetic mean/harmonic mean sit at opposite ends of a spectrum, that of \emph{$p$-means}.
		We give syntax and semantics for this quantitative predicate logic, and as example applications, we show how softmax is the quantitative semantics of argmax, and R\'enyi entropy/Hill numbers are additive/multiplicative semantics of the same formula.
		Indeed, the additive reals also fit into the story by exploiting the \emph{Napierian duality} $-\log \adj 1/\exp$, which highlights a formal distinction between `additive' and `multiplicative' quantities.
		Finally, we describe two attempts at a categorical semantics via enriched hyperdoctrines.
		We discuss why hyperdoctrines are in fact probably inadequate for this kind of logic.
	\end{abstract}

	\section{Introduction}
	In 1847, Boole \cite{boole_mathematical_1847} published his treatise on the algebra of logic, where he gave algebraic axioms to manipulate truth values as quantities.
	Ever since, mathematicians have been trying to go back a certain measure and formulate formal languages which unite the \emph{symbolic} and the \emph{quantitative} in a single system.
	On the one hand, such a system should enable clean symbolic reasoning, like for traditional formal languages; on the other, it should have the nuance of the continuously-valued real number.

	There is more to it, however, than just these two aspects: after all, mathematicians started investigating many-valued logic a long time ago, with the work of {\L}ukasiewicz \cite{lukasiewicz_o_1920} on the one hand, and Scott \cite{scott1967proof} on the other.
	Recently, applying that technology to the case of the reals, seen as a poset of truth values, Figueroa and van der Berg studied the topos of real-valued sets in \cite{Figueroa2022}.

	However, the mix of symbolic and quantitative that mathematicians \emph{intue} to be needed is different than just having continuum-many truth values.
	It also involves giving arithmetic operations, like $+$ and $\times$, a logical status.
	Probability theory gets very close to such a logic, with probability measures assigning belief degrees to propositions which are then summed or multiplied in the ways every mathematician is familiar with, and there are indeed axiomatizations of probability theory which brings to light their \emph{rational\footnotemark~content} \cite{cox_probability_1946}.%
	\footnotetext{As in: pertaining reasoning.}
	However, probability theory lacks a formal or syntactic aspect to it, so that the probabilist must reason \emph{with the theory} instead of \emph{within the theory}.

	Something different is the \emph{real-valued generalised logic} of Lawvere \cite{lawvere_metric_1973}, which tries to take seriously the idea that operations on the reals are logical connectives which can be manipulated symbolically while targeting a quantitative semantics.
	Lawvere's ideas have been a great source of inspiration in category-theoretic circles, though only very recently Bacci, Mardare, Panangaden and Plotkin picked up where Lawvere left and started developing his logic \cite{bacci_propositional_2023,bacci_polynomial_2024}.

	Parallel to these coherent attempts, we must mention the vast development on fuzzy logic and fuzzy set theory \cite{zadeh_fuzzy_1965,zadeh_fuzzy_1988,barr_variable_1985} as well as the growing literature on differential logics with continuous semantics \cite{diff1,diff2}.
	Also various works dealing with probability have foreshadowed links between the symbolic and the quantitative, such as Kozen's probabilistic PDL \cite{kozen_probabilistic_1983} where predicates are measurable functions and satisfaction is measured by integration (not unlike what we are going to do in this work).
	Other works in this space that we would like to mention are \cite{grandis_categories_2007,dagnino_logical_2022,dagnino_quantitative_2022}.

	This work started as an attempt to understand the semiotics\footnote{Meaning not just the crude mathematical definition but the set of latent meaning and informal intuitions that a mathematical idea bears and is communicated with.} of $\softmax$, which is used and sometimes advertised as a `quantitative' (alternatively: `smooth', `continuous', `fuzzy') version of $\argmax$.
	Trying to understand in what sense is this intuition true led us to analyze the structure of the extended positive reals and we discovered a logical structure of \emph{degrees of belief} which partially recovers known observations (such as \cite{bacci_polynomial_2024,shulman_affine_2022,grandis_categories_2007,kozen_probabilistic_1983}), but are also original to some extent.

	Nonetheless, this work is the proverbial tip of the iceberg.
	We only got glimpses of the logic we are about to describe, and this concerns syntax and semantics alike.
	Specifically, we lack a proof system, which limits the reach of the justifications we can give formally to symbolic manipulations.

	\subsection{Outline of the paper.}
	In \cref{sec:reals} we analyze the extended positive real numbers $\MulReals$ and single out the algebraic structure it naturally hosts, in the form of three generations of logical connectives---a nonlinear one and two linear ones (additive and multiplicative).
	We single out positive reals as a \emph{multiplicative} world, dual to the \emph{additive} world which is instead given by the full extended real line $\AddReals$.
	Taking seriously this duality is an important conceptual idea of this work.

	In the multiplicative reals, the linear additive connectives given by \emph{sum} and its de Morgan dual, \emph{harmonic sum} \cite{grandis_categories_2007}, form the basis for \cref{sec:quants}, where we notice two things.
	First, both kinds of sums are actually part of a spectrum, well known in real analysis, of $p$-sums, given by conjugating $+$ with the non-trivial automorphisms of $\MulReals$.
	Second, by normalizing these sums we get $p$-means, readily generalised to $p$-integrals over measure spaces (since means are integrals over probability spaces) and these behave like \textbf{soft quantifiers} for predicates bounded by a measure space, reducing to essential quantifiers (those minding only subsets of stricitly positive measure) for $p=\infty$.
	In  \cref{sec:wild} we show our main result, which is that softmax and argmax are semantics of the same sentence.

	In \cref{sec:cat.semantics} we try to frame integrals as quantifiers in the sense of Lawvere \cite{lawvere_adjointness_1969}, thus as adjoints to reindexing.
	We first attempt to do so na\"ively, and then by deploying the theory of enriched indexed categories of Shulman \cite{shulmanEnrichedIndexedCategories2013} to come up with an enriched hyperdoctrine of quantitative predicates.
	We fail both times in telling ways, which is why we describe our failure: idempotency at the meta-level, manifesting in different forms, ends up ruining our plans in both cases.


	\subsection{Acknowledgments.}
	I am deeply indebted to Ekaterina Komendatskaya, Bob Atkey, and Radu Mardare for their encouragement, support, and advice regarding this work.
	We thank Toby St Clere Smithe, Owen Lynch, David Jaz Myers, Kevin Carlson, Paolo Perrone, Frederik Nordvall-Forsberg and Neil Ghani for many conversations which helped shape this investigation.
	We thank Zoltan A. Kocsis and Maya M. for spotting some fairly conspicuous typos in earlier drafts of this work, Zanzi Mihejevs for suggesting the link with almost everywhere quantification, and Jairo Miguel Marulanda-Giraldo for helpful comments and corrections.

	\section{The extended positive real numbers}
	\label{sec:reals}

	The algebraic structure of the real numbers seems a fairly simple matter, yet it took us quite some time to reach the definition below.
	The reason is, as always in these cases, that many different structures are at play and one has to districate the various tendrils to reach a satisfying state of affairs.

	For the reals, there are two major confounding factors: the mix of additive and multiplicative structure and the role played by $0$ and $\pm\infty$.
	In high school, and throughout college, we are taught that multiplication and addition (or, more accurately, subtraction and division) do not play nice with these extremal values.
	Specifically, $0\infty$ could be both $0$ and $\infty$, and $\infty - \infty$ could be both $+\infty$ or $-\infty$.

	We make sense of this situation by first separating additive and multiplicative algebraic structures, realizing they are two faces of the same coin, related by \emph{Napierian duality} $-\log \adj 1/\exp: \AddReals \cong \MulReals$.
	Then, focusing on the multiplicative reals, we resolve the conflict of defining $0\infty$ by weakening the structure of $\MulReals$ from `compact closed' (thus admitting multiplicative inverses) to `$*$-autonomous' (admitting an inverse-like duality) \cite{barr_autonomous_1979}.
	In practice, this does not weaken the structure but perfectly captures the sense in which $\MulReals$ allows taking inverses.

	We start by expounding the latter structure:

	\begin{definition}
		The poset $(\PosReals, \leq)$ is a $\ast$-autonomous quantale, which we call \textbf{multiplicative reals} and denote by $\MulReals$, given by the following structure:
		\begin{enumerate}
			\item it has all small joins $\bigvee_i a_i$, given by suprema;
			\item it is equipped with a tensor product $\tensor : \PosReals \times \PosReals \to \PosReals$, defined by multiplication extended with the rules
			\begin{equation}
				\forall a \in (0,\infty],\ a\tensor\infty = \infty, \quad 0\tensor\infty = 0,
			\end{equation}
			and satisfying the properties
			\begin{enumerate}
				\item $\tensor$ distributes over joins: $k \tensor \bigvee_i a_i = \bigvee_i k \tensor a_i$,
				\item $\tensor$ is commutative, associative and unital, with unit $1 \in \PosReals$.
			\end{enumerate}
			\item it is equipped with a duality $(-)^* : \PosReals \to \PosReals^\op$, defined as
			\begin{equation}
				\forall a \in (0,\infty),\ a^* := 1/a, \quad 0^* = \infty,\quad \infty^*=0,
			\end{equation}
			satisfying the property
			\begin{equation}
				\forall a \in \PosReals,\ a \tensor b \leq c^* \iff a \leq (b \tensor c)^*.
			\end{equation}
		\end{enumerate}
	\end{definition}

	\begin{remark}
		The definitions of $a \tensor \infty$ and $0 \tensor \infty$ are dictated by the requirement $\tensor$ preserves joins.
	\end{remark}

	We remind the reader that the law $(a\tensor b)^* = a^*\tensor b^*$ is not required to hold.
	In fact, it does not hold in $\MulReals$, since $0 = 0\tensor\infty = \infty^*\tensor0^* \neq (\infty\tensor 0)^* = \infty$.
    However, it holds for all $0 < a,b < \infty$.

	The reason is that a $\ast$-autonomous category (see \cite{barr_autonomous_1979,blanco_bifibrations_2020}) like $\MulReals$ really features two distinct monoidal products, where the second is de Morgan dual to the first: $a \tensor^* b := (a^* \tensor b^*)^*$.
	In the case of $\MulReals$, $\tensor^*$ is again given by multiplication on all finite values but, due to the fact $\tensor^*$ is a monoidal product on the dual order $\MulReals^\op$, preservation of joins now forces the definition $0\tensor^* \infty = \infty$.

	The two structures lax-linearly distribute \cite{cockett_weakly_1997} over each other, meaning $a \tensor (b \cotensor c) \leq (a \tensor b) \cotensor c$, although equality holds for all finite values of $a,b,c$ and only breaks for e.g.~$(a,b,c)=(0,0,\infty)$.

	We summarize these two operations here:
	\begin{equation}
		\begin{tabular}{c|ccc}
			$a \tensor b$ & $0$ & $a \in (0,\infty)$ & $\infty$\\
			\cline{1-4}
			$0$ 			   & $0$ & $0$ 		& $0$\\
			$b \in (0,\infty)$ & $0$ & $ab$		& $\infty$\\
			$\infty$ 		   & $0$ & $\infty$ & $\infty$
		\end{tabular}
		\hspace*{10ex}
		\begin{tabular}{c|ccc}
			$a \cotensor b$ & $0$ & $a \in (0,\infty)$ & $\infty$\\
			\cline{1-4}
			$0$ 		 	   & $0$ 		& $0$ 	   & $\infty$\\
			$b \in (0,\infty)$ & $0$ 		& $ab$	   & $\infty$\\
			$\infty$ 		   & $\infty$ 	& $\infty$ & $\infty$
		\end{tabular}
	\end{equation}

	Notice, moreover, that exponentiation defines an action of $\MulReals$ on itself:
	\begin{equation}
	\label{eq:exponentiation}
		\begin{tabular}{c|ccccc}
			$k \cdot a$ & $0$ & $a \in (0,1)$ & $1$ & $a \in (1, \infty)$ & $\infty$\\
			\cline{1-6}
			$0$ & $1$ & $1$ & $1$ & $1$ & $1$\\
			$k \in (0, \infty)$ & $0$ & $a^k$ & 1 & $a^k$ & $\infty$\\
			$\infty$ & $0$ & $0$ & $1$ & $\infty$ & $\infty$
		\end{tabular}
	\end{equation}
	With this definition, $a^1 =a$ and $(a^h)^k = a^{h \tensor k}$ (in particular, $(a^0)^\infty = a^{0 \tensor \infty}$).
	We note that this action is self-dual, i.e.~$((a^*)^k)^* = a^k$, for finite values of $k$.

	Of crucial importance is the presence of \textbf{division}:
	\begin{equation}
		\begin{tabular}{c|ccc}
			$a \multimap b$ & $0$ & $a \in (0,\infty)$ & $\infty$\\
			\cline{1-4}
			$0$				   & $\infty$ & $0$ 	  & $0$\\
			$b \in (0,\infty)$ & $\infty$ & $b/a$	  & $0$\\
			$\infty$		   & $\infty$ & $\infty$  & $\infty$
		\end{tabular}
	\end{equation}
	This definition is dictated by the requirement $a \multimap -$ be right adjoint to $- \tensor a$.
	For $*$-autonomous quantales, the identities $a \multimap 1 = a^*$ and $a \multimap  b = a^* \cotensor b$ hold.
	We slightly abuse notation from now on, using fraction notation $b/a$ or $\dfrac{b}a$ to denote the above operation.

	Finally, on $\MulReals$ we can define a second join-preserving monoidal operation, \textbf{sum}, which is the obvious extension of addition as shown below left.
	Of great interest is also its de Morgan dual, \textbf{harmonic sum}, shown below right.
	\begin{equation}
		\begin{tabular}{c|ccc}
			$a \add b$ & $0$ & $a \in (0,\infty)$ & $\infty$\\
			\cline{1-4}
			$0$ 		 	   & $0$	  & $a$		 & $\infty$\\
			$b \in (0,\infty)$ & $b$	  & $a+b$	 & $\infty$\\
			$\infty$		   & $\infty$ & $\infty$ & $\infty$
		\end{tabular}
		\hspace*{10ex}
		\begin{tabular}{c|ccc}
			$a \coadd b$ & $0$ & $a \in (0,\infty)$ & $\infty$\\
			\cline{1-4}
			$0$ 		 	   & $0$ & $0$					  & $0$\\
			$b \in (0,\infty)$ & $0$ & $\dfrac{1}{1/a + 1/b}$ & $b$\\
			$\infty$		   & $0$ & $a$ 					  & $\infty$
		\end{tabular}
	\end{equation}
	Harmonic sum is associative, commutative, and unital with unit $\infty$.
	If $(a_i)_{i \in I} \in \MulReals$ is a finite set of numbers, we write $\hsum{i \in I} a_i$ for their harmonic sum.

	\subsection{Very linear logic.}
	Hence $\MulReals$ supports a very rich structure: not only a tensor $\tensor$ and a sum $\add$, but also joins $\lor$ and a duality $(-)^*$ which doubles them all.
	We summarise it in \cref{table:prop-mult-ll}.

	\begin{table}[htbp]
		\centering
		\begin{tabularx}{\textwidth}{|c|c|*{3}{>{\centering\arraybackslash}X|}}
			\hline
			\multicolumn{2}{|c|}{} & & \multicolumn{2}{c|}{\textbf{linear}}\\
			\cline{1-2}\cline{4-5}
			\multicolumn{2}{|c|}{\textbf{polarity}} & \multirow{-2}{*}{\textbf{nonlinear}} & \textbf{additive} & \textbf{multiplicative}\\
			\hline
			\multirow{2}{10ex}[-2ex]{\centering\parbox[c]{10ex}{\centering\textbf{duality}\\$a^* := 1/a$}} & \textbf{positive} & $\begin{aligned} \false &:= 0 \\[-1.5ex] a \lor b &:= \max\{a,b\} \end{aligned}$ & $\begin{aligned} \Zero &:= 0 \\[-1.5ex] a \add b &:= a + b \end{aligned}$ & $\begin{aligned} \One &:= 1 \\[-1.5ex] a \tensor b &:= ab,\\[-1.5ex] 0 \tensor \infty :&= 0 \end{aligned}$\\[2ex]
			\cline{2-5}
			& \textbf{negative} & $\begin{aligned} \true &:= \infty \\[-1.5ex] a \land b &:= \min\{a,b\} \end{aligned}$ & $\begin{aligned} \top &:= \infty \\[-1.5ex] a \with b &:= a +^* b \end{aligned}$ & $\begin{aligned} \bot &:= 1 \\[-1.5ex] a \cotensor b &:= ab,\\[-1.5ex] 0\cotensor \infty &:= \infty \end{aligned}$\\[2ex] 
			\hline
		\end{tabularx}

		\vspace*{2ex}
		\textbf{scalar modality} \hspace*{10ex} $k \cdot a := a^k, \qquad k \in \MulReals$

		\caption{Multiplicative semantics of very linear logic in $\MulReals$.}
		\label{table:prop-mult-ll}
	\end{table}

	The various monoidal products distribute over each other in distinct ways:
	\begin{eqalign*}
	\label{eq:distrib}
		\textbf{lax linear} && a \cotensor (b \tensor c) &\leq (a \cotensor b) \tensor c & \One &\leq \bot\\
		\textbf{nonlinear} && a \tensor (b \add c) &= (a \tensor b) \add (a \tensor c) & a \tensor 0 &= 0\\
		\textbf{nonlinear} && a \tensor (b \lor c) &= (a \tensor b) \lor (a \tensor c) & a \tensor \false &= \false\\
		&& a \add (b \lor c) &= (a \add b) \lor (a \add c) & a \add \true &= \true
	\end{eqalign*}

	If we interpret these as logical connectives, we are faced with a version of linear logic with \emph{three} generations of connectives: \textbf{nonlinear}, \textbf{linear additive} and \textbf{linear multiplicative}.
	Linear ones behave like the multiplicatives in classical linear logic, and the nonlinear ones like the additives.
	But also linear multiplicatives and linear additives interact in the same way they do in linear logic, hence their names.
	Like in classical linear logic, all connectives come in two polarities which are exchanged by a formal duality, which also functions as negation.

	For lack of a better name, we refer to this extended system as \textbf{very linear logic}.

	Unfortunately, it is not clear to us in which way the rules for linear additive connectives should distinguish them from linear multiplicative ones.
	In fact, we are not aware of other occurrences of this kind of (likely substructural) logic.

	Additionally, while we do not know whether it is legitimate to call this logic linear: we do not have a proof system that justifies calling $\add$ and $\tensor$ `linear connectives'.

	\subsection{Napierian duality.}
	Indeed, \cite{bacci_propositional_2023} introduced $\AL$, a propositional logic valued in $([0,\infty], {\geq}, \add)$, and recently \cite{bacci_polynomial_2024} introduced $\PL$ which extends $\AL$ with $\tensor$ (using our notation).
	The link with the structure we are contemplating is given by the following well-known diagram of quantales:
	\begin{equation}\label{eq:napier}
		\begin{tikzcd}[ampersand replacement=\&,row sep=scriptsize]
			\&[-4ex] {\textbf{multiplicative}} \&\& {\textbf{additive}} \\[-1ex]
			{\textbf{negative}} \& \NegMulReals \&\& {\NegAddReals^\op} \\[1ex]
			\& \MulReals \&\& {\AddReals^\op} \\[1ex]
			{\textbf{positive}} \& \PosMulReals \&\& {\PosAddReals^\op}
			\arrow[""{name=0, anchor=center, inner sep=0}, "{{-\log}}"', curve={height=6pt}, from=3-2, to=3-4]
			\arrow[""{name=1, anchor=center, inner sep=0}, "{{1/\exp}}"', curve={height=6pt}, from=3-4, to=3-2]
			\arrow[hook, from=4-4, to=3-4]
			\arrow[hook, from=4-2, to=3-2]
			\arrow[""{name=2, anchor=center, inner sep=0}, "{{-\log}}"', curve={height=6pt}, from=4-2, to=4-4]
			\arrow[""{name=3, anchor=center, inner sep=0}, "{{1/\exp}}"', curve={height=6pt}, from=4-4, to=4-2]
			\arrow[""{name=4, anchor=center, inner sep=0}, "{-\log}"', curve={height=6pt}, from=2-2, to=2-4]
			\arrow[""{name=5, anchor=center, inner sep=0}, "{1/\exp}"', curve={height=6pt}, from=2-4, to=2-2]
			\arrow[hook, from=2-2, to=3-2]
			\arrow[hook, from=2-4, to=3-4]
			\arrow["{1/-}", shift left=5, curve={height=-18pt}, dotted, tail reversed, from=4-2, to=2-2]
			\arrow["{-(-)}"', shift right=5, curve={height=24pt}, dotted, tail reversed, from=4-4, to=2-4]
			\arrow["{\begin{smallmatrix}\textbf{Napierian}\\\textbf{duality}\end{smallmatrix}}", tail reversed, from=1-2, to=1-4]
			\arrow["{\begin{smallmatrix}\textbf{de Morgan}\\\textbf{duality}\end{smallmatrix}}"', tail reversed, from=2-1, to=4-1]
			\arrow["\wr"{marking, allow upside down}, draw=none, from=1, to=0]
			\arrow["\wr"{marking, allow upside down}, draw=none, from=3, to=2]
			\arrow["\wr"{marking, allow upside down}, draw=none, from=5, to=4]
		\end{tikzcd}
	\end{equation}

	The logics $\AL/\PL$ deal with the positive additive reals, while the very linear logic above is interpreted in the multiplicative reals.
	The reason we favour working multiplicatively rather than additively is that the three generations of structure we expounded above are more explicit in this setting.
	Indeed, the linear additive structure on $\AddReals$ is awkward, arising because Napierian duality can be used to transport the structure from $\MulReals$ rather than in a natural manner---albeit this remains subjective.
	In \cref{table:prop-add-ll} we give the result of doing so: the linear multiplicative fragment of the additive reals is given by addition, so that the linear additive operation has to be something in-between, an operation known as log-sum-exp (LSE) \cite{nielsen_guaranteed_2016}, also \textbf{softplus} \cite{glorot_deep_2011}, which we prefer.

	\begin{table}[tbp]
		\centering
		\begin{tabularx}{\textwidth}{|c|c|>{\centering\arraybackslash}X|>{\centering\arraybackslash\hsize=23ex}X|>{\centering\arraybackslash}X|}
			\hline
			\multicolumn{2}{|c|}{} & & \multicolumn{2}{c|}{\textbf{linear}}\\
			\cline{1-2}\cline{4-5}
			\multicolumn{2}{|c|}{\textbf{polarity}} & \multirow{-2}{*}{\textbf{nonlinear}} & \textbf{additive} & \textbf{multiplicative}\\
			\hline
			\multirow{2}{10ex}[-2ex]{\centering\parbox[c]{10ex}{\centering\textbf{duality}\\$a^* := -a$}} & \textbf{positive} & $\begin{aligned} \false &:= -\infty \\[-1.5ex] a \lor b &:= \min\{a,b\} \end{aligned}$ & $\begin{aligned} \Zero &:= -\infty \\[-1.5ex] a \add b &:= -\log(\e^{-a} + \e^{-b}) \end{aligned}$ & $\begin{aligned} \One &:= 0 \\[-1.5ex] a \tensor b &:= a + b,\\[-1.5ex] -\infty \tensor \infty &:= +\infty \end{aligned}$\\[2ex]
			\cline{2-5}
			& \textbf{negative} & $\begin{aligned} \true &:= +\infty \\[-1.5ex] a \land b &:= \max\{a,b\} \end{aligned}$ & $\begin{aligned} \top &:= \infty \\[-1.5ex] a \with b &:= \log(\e^a + e^b) \end{aligned}$ & $\begin{aligned} \bot &:= 0 \\[-1.5ex] a \cotensor b &:= a +^* b,\\[-1.5ex] -\infty \cotensor \infty &:= -\infty \end{aligned}$\\[2ex]
			\hline
		\end{tabularx}

		\vspace*{2ex}
		\textbf{scalar modality} \hspace*{10ex} $k \cdot a := ka, \qquad k \in \MulReals$

		\caption{Additive semantics of very linear logic in $\AddReals$. The linear additive operations, given by $\add$ and $\coadd$, are called, respectively, \textbf{softplus} and \textbf{harmonic softplus}. Note the ordering is reversed compared to the usual order of the reals, though we agree with \cite{lawvere_metric_1973} and \cite{bacci_propositional_2023}.}
		\label{table:prop-add-ll}
	\end{table}

	Since we are ultimately interesting in understanding quantification, and $\sum$ is a clear candidate for the role of $\exists$, we are particularly interested in working in a setting where the linear additive structure is natural.
	Moreover, we prefer using both negative and positive values than limiting ourselves to positives as to have available not only degrees of truth but also of falsehood.

	\begin{remark}
	\label{rmk:two-worlds}
		The very existence of two `worlds' of quantitative operations, related by a perfect duality, is conceptually and mathematically stimulating.
		One can observe this duality manifesting in various occurrences, where multiplicative quantities (e.g.~probabilities, likelihood, perplexity, diversity, relative errors, ratios) and additive ones (e.g.~logits, information, entropy, energy, absolute errors, distance) interact in a controlled way.
		In particular this dichotomy makes one pay attention to the different meanings attached to real numbers and their operation.
		This is particularly intriguing in the twofold nature of sum: sum of additive quantities is different, conceptually, from the sum of multiplicative ones.
	\end{remark}

	\subsection{Degrees of truth.}\label{sec:sep}
	A few lines above we noted we can think of the numbers in $\MulReals$ as truth values.
	To convert these values back to classical ones, we use maps into the quantale of (meta-level) truth values $\Prop_\land$:

	\begin{definition}
		A \textbf{separator} on a quantale $Q_\tensor$ is a map of quantales $\lfloor-\rfloor:Q_\tensor^\op \to \Prop_\land$.
	\end{definition}

	A separator casts each $q \in Q_\tensor$ to an actual truth value $\lfloor q\rfloor \in \Prop_\land$.
	The contravariance is due to the fact such a map needs to preserve joins, which in $Q_\tensor^\op$ corresponds to meets in $Q_\tensor$.
	Thus a separator defines a subset $S = \{q \mid \lfloor q\rfloor = \bot\}$ of \textbf{qualitatively true} elements which is
	\begin{enumerate}
		\item \textbf{closed under arbitrary meets}, i.e. for all small families $(a_i)_{i \in I}$ of elements of $S$, $\bigwedge_i a_i \in S$,%
		\footnote{Since $\bigwedge \varnothing = \top$, this implies $S$ is \textbf{inhabited} too, i.e.$\top \in S$.}
		\item \textbf{upwards closed}, i.e. if $a \in S$ and $a \leq b$, then $b \in S$,
		\item \textbf{closed under $\tensor$}, i.e. if $a, b \in S$, then $a \tensor b \in S$,
	\end{enumerate}
	\emph{Vice versa}, any such subset defines a separator by the map $\_ \in S$, hence we refer as separators to both concepts.

	It is easy to classify the separators of $\MulReals$: by the second condition they must be intervals unbounded on the right, by the first they must be closed intervals, and the third condition rules out all the intervals of the form $[t,\infty]$ where $0 < t < 1$.
	Thus we are left with the \textbf{inconsistent separator} $[0,\infty]$, which considers everything to be true; and all the (principal) separators $[t, \infty]$ defined by a \emph{threshold for truth} $t \geq 1$ (possibly $t=\infty$).

	In this latter family, two separators have a distinguished status.
	The first is the \textbf{unitary} separator $U=[1,\infty]$.
	It corresponds to `positive numbers': indeed this (after passing through $-\log$) is the choice (implicitly) made for $\AL$ and $\PL$.
	In general quantales, it is the separator defined by the upperset of the unit of $\tensor$.

	The second is the \textbf{definite} separator $D=[\infty, \infty]$, and it corresponds to the attitude of considering any finite quantity of evidence insufficient to show truth.
	In general quantales, it is the separator defined by the upperset of the top element $\top$ (hence always a singleton).


	\section{Predicates and quantification}
	\label{sec:quants}

	As anticipated, sums emerge as natural candidates for quantifiers: they seem to aggregate evidence by summing up the values of a `quantitative predicate' $\varphi : I \to \MulReals$.
	This is hardly an original observation: the correspondence between the symbolic and the quantitative has long lived in the collective subconscious of mathematicians, but also surfaced consciously in places like e.g.~\cite{kozen_probabilistic_1983, perrone_kan_2021}.

	Let us highlight some easy-to-verify properties of $\psum[]{}$ and its dual, $\hsum{}$, which hint at their role as quantifiers:

	\begin{lemma}\label{lemma:har-sum}
		Harmonic sum
		\begin{enumerate}
			\item satisfies the identity:
			\begin{equation}
			\label{eq:inv-add-sum}
				\hsum{i\in I} \dfrac{a}{b_i} = \dfrac{a}{\psum[]{i \in I} b_i}
			\end{equation}
			\item satisfies the Fubini property:
			\begin{equation}
				\hsum{i \in I}\hsum{j \in J} a_{ij} = \hsum{(i,j) \in I \times J} a_{ij} = \hsum{j \in J}\hsum{i \in I} a_{ij},
			\end{equation}
			\item is homogeneous (i.e.~multiplication distributes over it):
			\begin{equation}
			\label{eq:homog}
				k \cotensor \hsum{i \in I} a_i = \hsum{i \in I} k \cotensor a_i, \quad k \tensor \hsum{i \in I} a_i = \hsum{i \in I} k \tensor a_i,
			\end{equation}
			implying moreover
			\begin{equation}
				\left(\hsum{i \in I} a_i\right) \tensor \left(\hsum{j \in J} b_j\right) = \hsum{i \in I}\hsum{j \in J} a_i \tensor b_j,
			\end{equation}
			\item is monotonic in the argument: if, for each $i \in I$, $a_i \leq b_i$, then
			\begin{equation}
				\hsum{i \in I} a_i \leq \hsum{i \in I} b_i,
			\end{equation}
			\item is antitonic in the index: when $J \subseteq I$, one has
			\begin{equation}
				\hsum{i \in I} a_i \leq \hsum{j \in J} a_{j}.
			\end{equation}
		\end{enumerate}
	\end{lemma}

	The identity \eqref{eq:inv-add-sum}, which we call the \textbf{fundamental property of harmonic sums}, can be written more suggestively as
	\begin{equation}
	\label{eq:fun-prop}
		\hsum{i} (b_i \multimap a) = \left(\psum[]{i} b_i\right) \multimap a,
	\end{equation}
	which is remindful of the usual $\forall i\, (b_i \to a) = (\exists i\, b_i) \to a$.
	Homogeneity can also be recast as
	\begin{equation}
	\label{eq:fun-prop2}
		\hsum{i} (b \multimap a_i)
		= b \multimap \left(\hsum{i} a_i\right),
	\end{equation}
	which is analogous to $\forall i\, (b \to a_i) = b \to (\forall i\, a_i)$.


	In fact, $\psum[]{}/\hsum{}$ and $\exists/\forall$ sit together in a spectrum:

	\begin{definition}[$p$-sum]
		For any $p \in (-\infty,\infty)$, $p \neq 0$, the \textbf{$p$-sum} of a finite set of numbers $(a_i)_{i \in I}$ is
		\begin{equation}
			\psum{i \in I} a_i := \left(\psum[]{i \in I} a_i^p\right)^{1/p}.
		\end{equation}
		To refer explicitly to the negative $p$ case, we talk of \textbf{harmonic $p$-sums} instead of $(-p)$-sums.
	\end{definition}

	\begin{remark}
	\label{rmk:hsums-as-conj}
		Like harmonic sum is the conjugate of arithmetic sum under taking reciprocals (i.e.~under the automorphism $(-)^*$ of $\MulReals$), $p$-sums are the conjugate of arithmetic sum under the scalar modalities $(-)^p$, thus making $p$-sums `modal' deformations of the usual sum.
	\end{remark}

	One extends the above definition to $p=\pm \infty$ (but not $p=0$) by taking suitable limits (see \cite[75]{mitrinovic1970analytic}).
	Then varying $p$ we recover sum ($p=1$), harmonic sum ($p=-1$), $\bigwedge$ ($p=-\infty$) and $\bigvee$ ($p=+\infty$), the latter two being the semantic counterpart of $\forall$ and $\exists$:
	\begin{equation}
		\begin{tikzcd}[sep=scriptsize, column sep=3ex]
			\bigwedge & \cdots &[-1.5ex] {\psum[-p]{}} & \cdots &[-1.5ex] {\hsum{}} &[1ex] {\textcolor{gray}{\not\exists}} &[1ex] {\psum[]{}} & \cdots &[-1.5ex] {\psum{}} & \cdots &[-1.5ex] \bigvee \\[-2ex]
			{-\infty} && {-p} && {-1} & 0 & 1 && p && \infty
			\arrow[no head, from=2-5, to=2-6]
			\arrow[no head, from=2-6, to=2-7]
			\arrow[no head, from=2-1, to=2-3]
			\arrow[no head, from=2-3, to=2-5]
			\arrow[from=2-9, to=2-11]
			\arrow[no head, from=2-9, to=2-7]
		\end{tikzcd}
	\end{equation}

	\begin{remark}
		We observe that all of the items of \cref{lemma:har-sum} still hold for general $p$-sums, with the exception of monotonicity in the indexing which reverses for positive $p$.
		In particular, the fundamental properties \eqref{eq:fun-prop} and \eqref{eq:fun-prop2} still put in relation $\psum[-p]{}$ with $\psum[p]{}$ for all $p \in [-\infty,\infty]$.
	\end{remark}

	For all $p \neq 0$, we have
	\begin{equation}
		\psum[-p]{} \leq \bigwedge \leq \bigvee \leq \psum[p]{},
	\end{equation}
	with the gap between linear additive and nonlinear narrowing as $p$ grows (see \cref{fig:charts}).

	Intuitively, compared to $\exists$, a $p$-sum looks up existence in a set by allowing `cumulative' effects to compensate.%
	\footnote{Similarly to how existential quantification in a topos is defined by `summing up' local truths}
	Dually, $\forall$ checks uniform satisfaction of a predicate while $\psum[-p]{}$ imposes extra `regularity' conditions.
	For instance, $\hsum{i \in I} a_i \geq 1$ implies that each $a_i \geq 1$, but that is not sufficient: they have to be sufficiently larger than $1$ as to also satisfy $\psum[]{i \in I} 1/a_i \leq 1$!
	Increasing $p$ makes the latter condition weaker, until we reach that of $\bigwedge$.

	\subsection{Justifying the means.} 
\begin{figure}
		\centering
		\begin{tikzpicture}
	\begin{axis}[
		width={.7\textwidth},
		height={.3\textheight},
		xlabel={$p$},
		ylabel={truth value},
		legend columns=1,
		legend cell align={left},
		legend style={
			at={(1.05,0.5)},
			anchor=west,
			draw=none,
			/tikz/every even column/.append style={column sep=0.5cm}
		},
		grid=both,
		xmin=1,
		xmax=30,
		ymin=0,
		ymax=3
	]

	\addplot[
		color=green!60!black,
		smooth,
		solid,
		line width=1pt
	] table [x=p, y=p_sum, col sep=tab] {sample.csv};
	\addlegendentry{$p$-sum}

	\addplot [
		color=violet,
		smooth,
		solid,
		line width=1pt
	] table [x=p, y=harmonic_p_sum, col sep=tab] {sample.csv};
	\addlegendentry{harmonic $p$-sum}

	\addplot[
		color=red,
		smooth,
		dash pattern=on 6pt off 1.5pt,
		line width=1pt
	] table [x=p, y=p_mean, col sep=tab] {sample.csv};
	\addlegendentry{$p$-mean}

	\addplot[
		color=blue,
		smooth,
		dash pattern=on 6pt off 1.5pt,
		line width=1pt
	] table [x=p, y=harmonic_p_mean, col sep=tab] {sample.csv};
	\addlegendentry{harmonic $p$-mean}

	\addplot[
		color=black,
		smooth,
		dotted,
		line width=1pt
	] table [x=p, y=max, col sep=tab] {sample.csv};
	\addlegendentry{max}

	\addplot[
		color=black,
		smooth,
		dotted,
		line width=1pt
	] table [x=p, y=min, col sep=tab] {sample.csv};
	\addlegendentry{min}

	\end{axis}
\end{tikzpicture}
		\caption{A comparison of $p$-sums (solid lines) and $p$-means (dashed lines) in their behaviour as $p$ changes. There is a clear tendency towards max/min (dotted) as $p \to \infty$, and the convergence of the $p$-means to the geometric one as $p \to 0$.}
		\label{fig:charts}
	\end{figure}
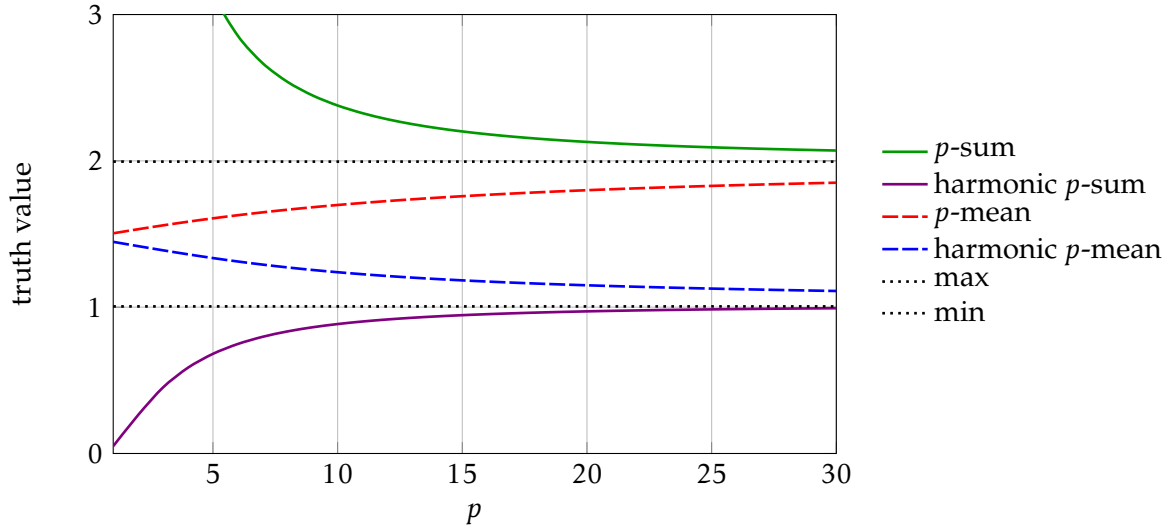

These characteristics are why $p$-sums are so useful, however, if left unbridled, they make $p$-sums depend strongly on the `size' of the quantification domain.
	Specifically, quantifying with a $p$-sum can yield arbitrarily high truth values depending only on how many `false' values we sum.
	This is apparent in \cref{fig:charts}, where you can see the $p$-sum going literally off the charts for small values of $p$.

	This is not necessarily a drawback, but there is an alternative worth exploring: we can compensate for size effects by taking \emph{means} instead of \emph{sums}.
	In fact, means are simply integrals over probability spaces:

	\begin{definition}
		For $p \in (-\infty,\infty)$ but $p\neq 0$, consider a probability space $(I, \de i)$ and a measurable function $a:I \to \PosReals$.
		We let the \textbf{$p$-mean} of $a$ be
		\begin{equation}
			\pmean[p]{i \in I} a(i)\, \de i := \left(\int_{i \in I} a(i)^p\, \de i\right)^{1/p}.
		\end{equation}
	\end{definition}

	Again, this definition can be extended to $p=\pm\infty$ by taking suitable limits.
	Unlike for sums, this time the limit exists even for $p=0$, and it is given on finite values by the \emph{geometric mean}.
	However, \textbf{when infinite values are involved, the geometric mean splits in two}, reflecting the two different multiplications on $\PosReals$:\footnotemark
	\footnotetext{
		The symbols $\pmean[+0]{}$ and $\pmean[-0]{}$ are formally defined by such limits, but can also be directly defined as multiplicative integrals \cite{bashirovMultiplicativeCalculusIts2008} by approximation by step functions.
		A step function is one of the form $\varphi=\sum_{k =0}^{N_\varphi} \varphi_k \tensor {\bf 1}_{E_k}$ where $N_\varphi < \infty$, the $E_k$ are a finite family of non-negligible measurable sets in $(I, \de i)$ and $\varphi_k \in (0,\infty]$; then $\pmean[-0]{i \in I} a(i)\,\de i = \sup_{\varphi \leq a} \bigotimes^{N_\varphi}_{k=0} \varphi_k^{\de i(E_k)}$, where the supremum ranges on all step functions bounded above by $a$.
		Replacing $\cotensor$ with $\tensor$ in the latter definition yields $\pmean[+0]{}$.}
	\begin{equation}
		\pmean[-p]{i \in I} a(i)\, \de i \conv[p \to 0^+] \pmean[-0]{i \in I} a(i)\, \de i, \qquad \pmean[p]{i \in I} a(i)\, \de i \conv[p \to 0^+] \pmean[+0]{i \in I} a(i)\,\de i.
	\end{equation}
	We call the first the (weighted) \textbf{conjunctive geometric mean}, and the second (weighted) \textbf{disjunctive geometric mean}.
	The difference can be grasped easily when e.g.~$I=2$ with the uniform measure and $a_1 = 0$, $a_2=\infty$: then $\pmean[+0]{i \in I} a(i)\,\de i = \infty$ while $\pmean[-0]{i \in I} a(i)\,\de i = 0$.
	Thus we can consider $\pmean[-0]{}$ and $\pmean[+0]{}$ infinitary versions of $\tensor$ and $\cotensor$, respectively, just like $\pmean[*]{}$ and $\pmean[]{}$ are infinitary versions of $\coadd$ and $\add$.

	We thus arrive at the following spectrum of means:
	\begin{equation}
		\begin{tikzcd}[sep=scriptsize, column sep=3ex]
			\bigwedge & \cdots &[-1.5ex] {\pmean[-p]{}} & \cdots &[-1.5ex] {\hmean{}} &[1ex] {\ \ \pmean[+0]{}\quad \pmean[-0]{}} &[1ex] {\pmean[]{}} & \cdots &[-1.5ex] {\pmean{}} & \cdots &[-1.5ex] \bigvee \\[-2ex]
			{-\infty} && {-p} && {-1} & 0 & 1 && p && \infty
			\arrow[no head, from=2-5, to=2-6]
			\arrow[no head, from=2-6, to=2-7]
			\arrow[no head, from=2-1, to=2-3]
			\arrow[no head, from=2-3, to=2-5]
			\arrow[from=2-9, to=2-11]
			\arrow[no head, from=2-9, to=2-7]
		\end{tikzcd}
	\end{equation}
	The extraordinary fact is that \textbf{$p$-means satisfy again the properties of \cref{lemma:har-sum}}, including Fubini, monotonicity in the argument and index, the fundamental relation and homogeneity.

	Unlike $p$-sums however, their values sit in between those of the traditional, nonlinear quantifiers, with the gap closing as $p$ grows (a trend very clear in \cref{fig:charts}):
	\begin{equation}
		\bigwedge \leq \pmean[-p]{} \leq \pmean[p]{} \leq \bigvee,
	\end{equation}
	Thus, unlike $p$-sums, \textbf{$p$-means compensate for size effects}.

	The intuition we propose is that sums and means are, respectively, unbounded and bounded quantifiers:
	\begin{eqalign}
		\exists i .\, a(i) &\analogous \psum[p]{i \in I} a(i)\\
		\exists (i \in I) .\, a(i) &\analogous \pmean[p]{i \in I} a(i)\,\de i \equiv \left(\psum[p]{i \in I} \dfrac1{|I|} \cdot a(i) \right)^{1/p}
	\end{eqalign}
	The mass element $\de i$ behaves as the bound $i \in I$ in a quantification.
	The fact that this is part of the quantifier can be seen in two different ways.
	First, on the syntactic side, we observe that $\de i$ does not participate in the duality relating $\psum[p]{i \in I}$ and $\psum[-p]{i \in I}$, like $i \in I$ does not for $\exists(i \in I)/\forall(i \in I)$:
	\begin{equation}
		\exists (i \in I) .\, a(i) = \neg \forall(i \in I) \neg a(i) \analogous \pmean[p]{i \in I} a(i)\,\de i = \left(\pmean[-p]{i \in I} a(i)^*\,\de i\right)^*.
	\end{equation}
	Secondly, on the semantic side, integrating with respect to $\de i$ is a more general operation than summing expressions such as $a(i) \cdot \de i$ since $\de i$ is a measure not a density (despite our suggestive notation), witnessing the fact that, generally speaking, it is part of the quantification itself and not of the argument.

	To sum it up, \textbf{we propose real-valued quantification should be obtained through integration}.
	The modalities $(\_)^p$ `soften out' integrals to give the necessary flexibility for applications. For $p \to \infty$, we recover `essential quantifiers':
	\begin{equation}
	\label{eq:ae.quant}
		\pmean[\infty]{i \in I} a(i)\,\de i = \esssup_{i \in I} a(i), \qquad \pmean[-\infty]{i \in I} a(i)\,\de i = \essinf_{i \in I} a(i).
	\end{equation}
	In particular, for a binary predicate $a : I \to \{0,1\}$, the first yields $1$ if $a$ is $1$ on a subset of strictly positive measure, and likewise the second gives $1$ if $a$ is $1$ almost everywhere.

	As known from real analysis, as $p$ increases, $p$-means give more and more importance to `oscillations', i.e. relative peaks and throughs---we say it gets a more \emph{vertical} character.
	Conversely, for small values of $p$ the norm has a more \emph{horizontal} character, downplaying the impact of oscillations.
	From a logical standpoint, we see both $\pmean[p]{}$ and $\pmean[-p]{}$ build up their value by adding up both `positive' ($>1$) and `negative' ($<1$) evidence, but as $p$ increases, (relatively) positive evidence is given more and more importance while (relatively) negative evidence is increasingly ignored.

	\subsection{A sketch of syntax and semantics.} To formalise the analogies sketched above, we need first to formalise what predicates are.

	We start by saying that each measure space $(I,\de i)$ (which we denote simply as $I$) is a \textbf{context} (or \textbf{domain of discourse}).
	The \textbf{multiplicative quantitative predicates} over each context $I$ are defined by the following inductive definition (given simultaneously for all contexts):
	\begin{eqalign}
	\label{eq:mult-pred}
		\LT(I)_\tensor \ni \varphi(i) :=\phantom{\mid\,}& \const \,\mid\, a\\
					 &\mid\, \varphi(i) \mathbin{\text{op}} \varphi(i)
					 \,\mid\, \varphi(i) \multimap \varphi(i)
					 \,\mid\, \varphi(i)^*\\
					 &\mid\, k \cdot \varphi(i)\\
					 &\mid\, \exists^p (j \in J).\, \psi(i,j)
					 \,\mid\, \forall^p (j \in J).\, \psi(i,j)
	\end{eqalign}
	where
	\begin{enumerate}[label=(\roman*)]
		\item $\const$ stands for any of the constants $\false, \true, \Zero, \One, \top, \bot$ introduced in \cref{table:prop-mult-ll},
		\item $a$ ranges in measurable functions $I \to \MulReals$, $\text{op} \in \{\lor, \land, \add, \coadd, \tensor, \cotensor\}$,
		\item $k \in [0,\infty)$,
		\item $p \in [0,\infty]$, and $\psi \in \LT_\tensor(I \times J)$, for $J$ another context.
	\end{enumerate}
	Note this is not a minimal grammar, i.e. different sentences correspond to the same function (since e.g.~constants are particular kinds of measurable functions, operators are interdefinable, etc.)
	Indeed, ideally, the sets $\LT(I)_\tensor$ would be further quotiented by syntactic rewrites or, even better, equipped with an entailment relation but it is not yet clear how this entailment would look like.
	We discuss this in \cref{sec:cat.semantics}.

	Still, the \textbf{multiplicative semantics} bracket $\sem{-}_\tensor$ can be defined inductively on $\LT(I)_\tensor$.
	It takes values in measurable functions $I \to \MulReals$, and since every measurable function determines a predicate too we sometimes terminologically conflate the two notions.

	The propositional fragment is covered by \cref{table:prop-mult-ll} where we gave semantics for constants, negation (the duality $(-)^*$), and binary operators (and $\sem{a}_\tensor = a$).
	The $p$-soft quantifiers are interpreted as
	\begin{eqalign}
	\label{eq:mult-quant}
		\sem{\exists^p(j \in J).\, \varphi(i,j)}_\tensor &:= \pmean[p]{j \in J} \sem{\varphi(i,j)}_\tensor\, \de j,\\
		\sem{\forall^p(j \in J).\, \varphi(i,j)}_\tensor &:= \phmean[-p]{j \in J} \sem{\varphi(i,j)}_\tensor\, \de j.
	\end{eqalign}
	As noted above, in the limit $p=\pm\infty$ the right hand side \emph{means} $\bigvee$ and $\bigwedge$, while for $p=0$ we still get distinct $\sem{\forall^0}=\pmean[-0]{}$ and $\sem{\exists^0} = \pmean[+0]{}$, the (weighted) conjunctive and disjunctive geometric means.


	We can define \textbf{additive quantitative predicates} by the same inductive definition as in \eqref{eq:mult-pred}, except now constants range in $\AddReals$.
	This defines a set $\LT(I)_\add$.
	There is a syntactic translation $-\log:\LT(I)_\add \to \LT(I)_\tensor$ given by applying \emph{the function} $-\log$ on constants and defined to commute with all the other formers (operations and quantifiers alike).
	The inverse syntactic translation $1/\exp:\LT(I)_\tensor \to \LT(I)_\add$ is defined analogously.

	Then \textbf{additive semantics} $\sem{-}_\add$ takes values in measurable functions $I \to \AddReals$.
	It can be defined from the multiplicative one by Napierian duality:
	\begin{equation}
	\label{eq:add-semantics}
		\sem{\varphi(i)}_\add := -\log \sem{1/\exp(\varphi(i))}_\tensor
	\end{equation}
	On the propositional fragment, the effect is described by \cref{table:prop-add-ll} and $p$-soft quantifiers become:
	\begin{eqalign}
	\label{eq:add-quant}
		\sem{\exists^p(i \in I).\, \varphi(i)}_\add &:= -\frac1p \log \int_{i \in I} \e^{-p\sem{\varphi(i)}_\add}\, \de i,\\
		\sem{\forall^p(i \in I).\, \varphi(i)}_\add &:= +\frac1p \log \int_{i \in I} \e^{+p\sem{\varphi(i)}_\add}\, \de i.
	\end{eqalign}

	Again, in the limit $p=\pm\infty$ these recover $\bigvee$ and $\bigwedge$, and for $p=0$ $\sem{\forall^0}$ is the (weighted) conjunctive arithmetic mean (defined from $\tensor$ in \cref{table:prop-add-ll}) while $\sem{\exists^0}$, the (weighted) disjunctive arithmetic mean (defined from $\cotensor$).

	As far as we know, these expressions have not been considered \emph{as quantifiers} before, though, as we are going to see now, they are not new to mathematics.

	\subsection{Quantifiers in the wild.}
	\label{sec:wild}

	We have one important reason to believe the quantifiers of \eqref{eq:add-quant} (and by extension the original multiplicative ones of \eqref{eq:mult-quant}) are interesting, which is their appearance in statistics and machine learning.

	\subsubsection{Argmax and softmax}
	The $\softmax$ operator is used in machine learning to turn a vector of real numbers into a probability distribution, and in statistical mechanics, where it is known as \emph{Gibbs} or \emph{Boltzmann distribution} to get a probability distribution out of a positive energy functional.

	\begin{definition}
		Let $f : X \to \MulReals$ be a function.
		Its \textbf{softmax} is the function $X \to \NegMulReals$
		\begin{equation}
		\label{eq:softmax}
			(\softmax f)(x^*) = \dfrac{f(x^*)}{\pmean[]{x \in X} f(x)\,\de x}.
		\end{equation}
	\end{definition}

	On the other hand, in the Boolean world we have a similar construction:
	\begin{equation}
		(\argmax f)(x^*) = \forall x \in X,\ f(x) \leq f(x^*).
	\end{equation}

	We get a `soft version' by replacing hard quantifiers with soft ones.
	Define:
	\begin{equation}
		(\psoftmax f)(x^*) = \sem{\forall^p (x \in X).\, f(x) \multimap f(x^*)}_\tensor = \pmean[-p]{x \in X} \dfrac{f(x^*)}{f(x)}\,\de x = \dfrac{f(x^*)}{\pmean[p]{x \in X} f(x)\,\de x}
	\end{equation}
	Then it is easy to see that $\softmax = \psoftmax[1]$, whereas $\argmax = \lfloor\psoftmax[\infty]\rfloor$:
	\begin{equation}
		\lfloor (\psoftmax[\infty] f)(x^*) \rfloor = \lfloor \bigwedge_{x \in X} f(x) \multimap f(x^*) \rfloor = \lfloor \bigwedge_{x \in X} f(x^*)/f(x) \rfloor = (\argmax f)(x^*).
	\end{equation}
	where the map $\lfloor-\rfloor:\MulReals \to \Prop_\land$ is the unitary separator of \cref{sec:sep}.
	Indeed, the ratio $f(x^*)/f(x)$ is greater or equal than $1$ iff $x^*$ is a maximum point, thus yielding $\top$ for those points. Else, it is alway strictly below $1$ and yields $\bot$.



	More commonly, $\softmax$ is applied to functions of arbitrary sign by first passing them through Napierian duality, i.e.~usually $f=\e^{-u}$ in \eqref{eq:softmax}.
	The decaying exponential is antitone which means that $\softmax$, for functions so constructed, puts more mass on small values of $u$ than large ones, making it more like a `softmin'.
	In any case, we still want to `extract back' the distribution obtained by $\softmax$, thus it stands to reason to apply $-\log$ to the result:
	\begin{equation}
		L(x^*) = -\log \softmax \e^{-u(x^*)} = - u(x^*) + \log \pmean[]{x \in X} \e^{-u(x)} \,\de x = \sem{\forall^1 (x \in X).\, u(x) \multimap u(x^*)}_\add.
	\end{equation}
	The quantity $L$ is called \textbf{log-likelihood}, and we just reconstructed it by purely formal considerations.


	\subsubsection{Entropy} Consider a predicate $\varphi:I \to \MulReals$.
	For simplicity, assume $\varphi$ is negative and unitary, meaning it lands in $\NegMulReals$ and $\int_i \varphi(i) = 1$.
	The \textbf{R\'enyi entropy} of order $p \in [0,\infty]$ \cite[Definition~4.3.1]{leinster_entropy_2021} is the quantity
	\begin{equation}
	\label{eq:renyi}
		H_p(\varphi) = \dfrac{1}{1-p} \log \int_{i \in I} \varphi(i)^p\,\de i.
	\end{equation}
	When $p=0,1,\infty$ \eqref{eq:renyi} must be interpreted in the limit.
	Notably, for $p=1$, $H_p(\varphi)$ recovers \textbf{Shannon entropy}.

	The quantity $H_p(\varphi)$ is the additive semantics of $p$-universally quantifying the negation of $-\log \varphi$, up to a multiplicative constant:
	\begin{equation}
	\label{eq:h}
		H_p(\varphi) = \sem{{\textstyle\frac{p}{1-p}} \cdot \left(\forall^p (i \in I). -\log \varphi(i)^*\right)}_\add
	\end{equation}
	This expression showcases some of the ideas regarding additive/multiplicative type discipline we talked about in \cref{rmk:two-worlds}.
	Indeed, the logarithm appearing in \eqref{eq:h} is very natural and can formally be justified as the necessary casting from the multiplicative world of $\varphi$ to the additive world---indeed, entropy is an additive quantity.

	In fact, one could define a purely multiplicative quantity:
	\begin{equation}
		D_p(\varphi) = \sem{{\textstyle\frac{p}{1-p}} \cdot \left(\forall^p (i \in I). \varphi(i)^*\right)}_\tensor.
	\end{equation}
	This is known as \textbf{$p$-diversity} (or \textbf{$p$-th Hill number}) \cite[Definition~4.3.4]{leinster_entropy_2021}.
	We deem it remarkable the fact this corresponds to a very simple logical formula.


	\section{Towards a categorical semantics}
	\label{sec:cat.semantics}
	The definition of quantifiers we gave above seems to be directly in the groove of an algebraic account of first-order predicate logic, where contexts are provided by measure spaces and algebras of predicates are given by the (wannabe) quantales $\LT(I)_\tensor$ and take semantics in analogous quantale of real-valued random variables.

	We focus now on the latter, since that is a setting we can test our guesses on.
	Define the $\MulReals$-enriched \emph{quantale of functionals} $\M^p(I, \MulReals)$ with entailment relation tentatively given by
	\begin{equation}
		\varphi \entails_I \psi := \sem{\forall^p (i \in I).\, \varphi(i) \multimap \psi(i)}_\tensor = \left(\int_{i \in I} \dfrac{\varphi(i)^p}{\psi(i)^p}\,\de i\right)^{-1/p}.
	\end{equation}
	This seems to define the action on objects as a functor $\M^p:\Meas^\op \to \VCat{\MulReals}$. The action on measure non-increasing maps $f:J \to I$ is readily given by precomposition, which is indeed an enriched functor $f^*:\M^p(I) \to \M^p(J)$, since, for every $\varphi,\psi \in \M^p(I)$, we have
	\begin{equation}
		\varphi \entails_I \psi \leq f^*\varphi \entails_J f^*\psi \quad\rightsquigarrow\quad \left(\int_{i \in I} \dfrac{\varphi(i)^p}{\psi(i)^p}\,\de i\right)^{-1/p} \leq\ \left(\int_{j \in J} \dfrac{\varphi(f(j))^p}{\psi(f(j))^p}\,\de j\right)^{-1/p}.
	\end{equation}
	The proof hinges on the fact that, if $u:I \to \MulReals$ is a measurable function, then $\int_j u(f(j))\,\de j = \int_i u(i)\, \de (f_* \de j) \leq \int_i u(i)\, \de i$ (and the inequality reverses when we take reciprocals, as we do above).

	Moreover, we see enriched left and right adjoints to reindexing are quantification, a property that famously characterises quantifiers \cite{lawvere_adjointness_1969}.
	For simplicity, let $f=\pi_I :I \times K \to I$, then we have $\pexists{k} \adj \pi_I^*$ since for all $\rho \in \M^p(I \times K)$, $\psi \in \M^p(I)$, we have
	\begin{equation}\label{eq:adj-quants}
		\pmean[p]{k \in K} \rho(-,k)\,\de k \entails_I \psi = \rho \entails_{I \times K} \pi_I^*\psi
	\end{equation}
	which unpacks to
	\begin{equation}
		\phmean{i \in I} \left(\pmean{k \in K} \rho(i,k)\,\de k\right) \multimap \psi(i)\,\de i
		=
		\phmean{\begin{smallmatrix}i \in I\\k \in K\end{smallmatrix}} \rho(i,k) \multimap \psi(i)\,\de k\,\de i
	\end{equation}
	and this identity is true by the fundamental property of harmonic $p$-means.
	Similarly, we can conclude $\pi_I^* \adj \pforall{k}$.
	Thus it would seem that $\M^p$ could be used a semantics for first-order very linear logic with $\exists^p/\forall^p$ sent to $\pexists{}/\pforall{}$.

	Unfortunately, this approach is broken: the entailment relation we proposed is not a proper enriched preorder since it does not satisfy reflexivity nor transitivity.
	Reflexivity says
	\begin{equation}
		\varphi \entails_I \varphi = \left(\int_{i \in I} \dfrac{\varphi(i)^p}{\varphi(i)^p}\,\de i\right)^{-1/p} =\ |I|^{-1/p} \geq 1
	\end{equation}
	and it is satisfied only when $|I| \leq 1$, thus when $I$ is a probability spaces.
	But even if restricted to those, transitivity still breaks. Fixing $p=1$ for simplicity, notice first that
	\begin{equation}
		(\varphi \entails_I \psi) \tensor (\psi \entails_I \sigma) \leq (\varphi \entails_I \sigma) \iff (\varphi \entails_I \psi)^{-1} \tensor (\psi \entails_I \sigma)^{-1} \geq (\varphi \entails_I \sigma)^{-1}.
	\end{equation}
	However, unpacking definitions, we see that:
	\begin{equation}
		\int_{i \in I} \int_{i' \in I} \dfrac{\varphi(i)}{\psi(i)}\dfrac{\psi(i')}{\sigma(i')}\,\de i\, \de i' \geq \int_{i \in \Delta_I} \dfrac{\varphi(i)}{\cancel{\psi(i)}}\dfrac{\cancel{\psi(i)}}{\sigma(i)}\,{\de i\,\de i} \ngeq \int_{i \in I} \dfrac{\varphi(i)}{\sigma(i)}\,{\de i}.
	\end{equation}
	From a logical standpoint, we are witnessing the failure of $\entails_I$ to satisfy the identity axiom and the cut rule, at least if expressed na\"ively.


	\subsection{The indexed way.}
	One might try to fix the previous attempt by varying the base of enrichment along with the indexing---the rationale being that prematurely taking integrals when defining entailment relation is what breaks the approach.

	This kind of variable enrichment has been invented by Shulman in \cite{shulmanEnrichedIndexedCategories2013}, where he describes enrichment of an indexed category in an indexed monoidal category.

	In our case the basis of enrichment is the following (we stick to $p=1$):

	\begin{definition}
		For each $\sigma$-finite\footnote{Meaning it can be covered by a countable family of finite measure spaces: it is the technical requirement for the existence of Radon--Nikodym derivatives which we employ below.} measure space $(I, \de i)$, the \textbf{quantale of functionals} on $I$ is the set $\M^1(I)$ of measurable functions $I \to \MulReals$ with essential pointwise ordering, meaning $\varphi \leq \psi$ iff $\varphi$ bounds $\psi$ from below pointwise almost-everywhere, and pointwise multiplication (using $\tensor$).
	\end{definition}

	Evidently, each quantale $\M(I)$ is closed and $*$-autonomous, and in fact has the same structure as $\MulReals$, inherited pointwise.
	The fact that $\M(I)$ is complete and cocomplete is \cite[Lemma~2.6.1]{meyer-nieberg_banach_1991}.
	Notice, in $\M(I)$, $\varphi \iso \psi$ when $\varphi$ equals $\psi$ almost everywhere on $I$.

	We would like to define the \textbf{indexed quantale of functionals} as the pseudofunctor
	\begin{equation}
		\M : \Meas^\op \longto \Quant
	\end{equation}
	defined as above on objects and acting on measure non-increasing functions $f:I \to J$ by precomposition. Thus $f^* : \M(J) \to \M(I)$ sends $\psi$ to $\psi(f)$.

	This definition is clearly well-posed, since $f^*$ is strict functorial as well as strict monoidal:
	\begin{equation}
		((f^*\varphi) \tensor (f^*\psi))(i) = \varphi(f(i)) \tensor \psi(f(i)) = f^*(\varphi \tensor \psi)(i).
	\end{equation}

	Now, still following \cite[Definition~4.1]{shulmanEnrichedIndexedCategories2013}, we self-enrich $\M$ to make it into an indexed $\M$-category.
	The data of such an enrichment, which we denote $\M_\M$, consists of
	\begin{enumerate}
		\item An $\M(I)$-enriched category $\M(I)_{\M(I)}$ for each $I \in \Meas$, and this enrichment is given by the fact $\M(I)$ is monoidal closed:
		\begin{equation}
			[\varphi, \psi]_{\M(I)_{\M(I)}} := \lambda (i \in I).\, \varphi(i) \multimap \psi(i).
		\end{equation}
		\item For each $f:I \to J$, a fully faithful $\M(I)$-functor $\overline{f^*} : \M(J)_{f^*} \to \M(I)_{\M(I)}$, where the domain is change of enrichment for $\M(J)$ along $f^*$.
		This means the hom-objects of this $\M(I)$-category are
		\begin{equation}
			[\varphi, \psi]_{\M(J)_{f^*}} := \lambda (i \in I).\, \varphi(f(i)) \multimap \psi(f(i)).
		\end{equation}
		Then $\overline{f^*}$ is defined on objects by precomposition with $f$, and on hom-objects is given by the identity
		\begin{equation}
			[\varphi, \psi]_{\M(J)_{f^*}} = [f^*\varphi, f^*\psi]_{\M(I)_{\M(I)}}.
		\end{equation}
	\end{enumerate}
	This data is also subject to coherence conditions which mostly trivialise for $\M_\M$, and thus we have a well-defined indexed quantale.

	Instead, one could try defining a covariant indexing
	\begin{equation}
		\M : \Meas \longto \Quant
	\end{equation}
	defined as above on objects and on measure non-increasing functions $f:I \to J$ defined as the monoidal monotone function $f_! : \M(I) \to \M(J)$ given by Radon--Nikodym derivative of the pushforward of measures:
	\begin{equation}
		f_!\varphi := \dfrac{\de f_*(\varphi \cdot \de i)}{\de j}
	\end{equation}
	where $\varphi \cdot \de i$ denotes the measure obtained by integrating $\varphi$ against $\de i$.

	Thus $f_!\varphi$ is the density of the function $\int_{i \in f^{-1}(-)} \varphi(i)\,\de i$.

	Pseudofunctoriality of $\M$ follows from applying the definition of Radon--Nikodym derivative and functoriality of pushforward of measures:
	\begin{equation}
		g_! f_! \varphi = \dfrac{\de g_*\left(\frac{\de f_*(\varphi \cdot \de i)}{\de j} \cdot \de j\right)}{\de k} = \dfrac{\de g_*(f_*(\varphi \cdot \de i))}{\de k} = \dfrac{\de (fg)_*(\varphi \cdot \de i)}{\de k} = (fg)_! \varphi.
	\end{equation}
	and the lack of strictness comes from the definition of Radon--Nikodym derivative only up to almost-everywhere equality.

	Lax monoidality is where the definition fails, in fact we have:
	\begin{eqalign}
		(f_! \varphi) \tensor (f_! \psi) \cdot \de j &= (f_! \varphi) \cdot f_*(\psi \cdot \de i)\\
		&= \int_{j \in (-)} \dfrac{\de f_*(\varphi \cdot \de i)}{\de j}(j)\, f_*(\psi \cdot \de i)\\
		&= \int_{i \in f^{-1}(-)} \dfrac{\de f_*(\varphi \cdot \de i)}{\de j}(f(i)) \, \psi(i)\,\de i\\
		&\textcolor{red}{\not\leq} \int_{i \in f^{-1}(-)} \varphi(i) \,\psi(i)\, \de i\\
		&= f_!(\varphi \tensor \psi) \cdot \de j
	\end{eqalign}
	We would like to use the inequality $\frac{\de f_*(\varphi \cdot \de i)}{\de j}(f(\bar i)) \leq \varphi(\bar i)$, i.e.~that $\int_{i \in f^{-1}(f(\bar i))} \varphi(i)\,\de i \leq \varphi(\bar i)$.
	But the fact $f^{-1}f(\bar i) \supseteq \{\bar i\}$ ca not be used to prove either direction of such an inequality since $\int_{i \in \{\bar i\}} \varphi(i) \,\de i \neq \varphi(i)$ in general.

	Thus the desired pseudofunctor $\M : \Meas \to \Quant$ does not exist.

	\begin{remark}
		One might object that $f_!$ is supposed to be left adjoint to $f^*$, thus colax monoidal, and indeed the fact we are comparing integrals over $f^{-1}f(\bar i) \supseteq \{\bar i\}$ seems to suggest this is the case.
		However, that direction also fails since we ca not compare evaluation at a point with integrals on a neighborhood in general.
	\end{remark}

	\section{Conclusions}
	\label{sec:conclusions}
	We described the syntax and semantics of a \emph{quantitative predicate logic} ($\QPL$), which naturally arises from the algebraic structures of the reals.
	We made our case that a logical interpretation can be given to some objects, like $\softmax$ and entropy, which are widely employed in non-symbolic quantitative reasoning for statistics, machine learning, physics, etc.

	We exhibited a syntax and a semantics for $\QPL$, but we could not come up with a convincing `proof theory' for it.
	In fact, the usual Boole--Gentzen--Lawvere algebraic framework for sequent calculus turned out to be inadequate, as testified by the failure of hyperdoctrines to capture $\QPL$.
	We are in need of a different algebraic structure corresponding to the rules of $\QPL$, capable of accomodating the lack of idempotency in its interpretation of identity and cut rules for entailment.

	Indeed, the influence a metatheory exerts on the theories it hosts is a fundamental aspect which is often underplayed.
	Hosting an unapologetically quantitative theory within a traditional qualitative metatheory seems very problematic, and we are thus left with a hard bootstrapping problem: that of formulating both a theory and a metatheory apt to host it at the same time.

	This obstacle is the only thing preventing us from declaring $\QPL$ a theory of quantitative \emph{reasoning}.
	As the present work stands, it is just a theory for quantitative \emph{syntax}.

	Finally, we have not covered an interesting use we can make of $\QPL$, which is to do category theory `enriched' in it.
	In that setting, one could formulate universal properties as equations some quantitative object has to satisfy.

	We had a taste of this idea in \eqref{eq:fun-prop}, the fundamental property of harmonic sum: the universal characterization of a quantitative object becomes an \emph{equation} that can be \emph{solved} to yield the desired object.%
	\footnote{In truth, a similar phenomenon holds in traditional settings, e.g.~a product is made of tuples because it is universally equipped with projections.}
	This enables to characterise some constructions: for instance, $\argmax$ can be characterised as a certain right Kan lift in the bicategory $\Rel$.
	The same universal property, written in $\QPL$, yields $\softmax$ as its unique solution.

	We are particularly interested in clarifying the universal properties of conditionals \cite{fritz_synthetic_2020} and Bayesian update \cite{jacobs_mathematics_2019,di_lavore_evidential_2023} in statistics: can we frame them as universal inductive reasoning rules, in a logical sense?

	\printbibliography

@misc{bacci_propositional_2023,
	title = {Propositional {Logics} for the {Lawvere} {Quantale}},
	url = {http://arxiv.org/abs/2302.01224},
	doi = {10.48550/arXiv.2302.01224},
	urldate = {2023-03-01},
	publisher = {arXiv},
	author = {Bacci, Giorgio and Mardare, Radu and Panangaden, Prakash and Plotkin, Gordon},
	month = feb,
	year = {2023},
	note = {arXiv:2302.01224 [cs]},
}

@misc{bacci_polynomial_2024,
	title = {Polynomial {Lawvere} {Logic}},
	url = {http://arxiv.org/abs/2402.03543},
	language = {en},
	urldate = {2024-02-12},
	publisher = {arXiv},
	author = {Bacci, Giorgio and Mardare, Radu and Panangaden, Prakash and Plotkin, Gordon},
	month = feb,
	year = {2024},
	note = {arXiv:2402.03543 [cs]},
}

@book{mitrinovic1970analytic,
  title={Analytic inequalities},
  author={Mitrinovic, Dragoslav S and Vasic, Petar M},
  volume={61},
  year={1970},
  publisher={Springer}
}

@article{cockett_weakly_1997,
	title = {Weakly distributive categories},
	volume = {114},
	doi = {10.1016/0022-4049(95)00160-3},
	journal = {Journal of Pure and Applied Algebra},
	author = {Cockett, Robin and Seely, R.A.G.},
	month = nov,
	year = {1997},
	pages = {133--173},
}

@article{blanco_bifibrations_2020,
	series = {The 36th {Mathematical} {Foundations} of {Programming} {Semantics} {Conference}, 2020},
	title = {Bifibrations of {Polycategories} and {Classical} {Linear} {Logic}},
	volume = {352},
	issn = {1571-0661},
	url = {https://www.sciencedirect.com/science/article/pii/S1571066120300499},
	doi = {10.1016/j.entcs.2020.09.003},
	urldate = {2024-03-20},
	journal = {Electronic Notes in Theoretical Computer Science},
	author = {Blanco, Nicolas and Zeilberger, Noam},
	month = oct,
	year = {2020},
	pages = {29--52},
}

@book{barr_autonomous_1979,
	address = {Berlin, Heidelberg},
	series = {Lecture {Notes} in {Mathematics}},
	title = {*-{Autonomous} {Categories}},
	volume = {752},
	isbn = {978-3-540-09563-7 978-3-540-34850-4},
	url = {http://link.springer.com/10.1007/BFb0064579},
	urldate = {2024-03-20},
	publisher = {Springer},
	author = {Barr, Michael},
	year = {1979},
	doi = {10.1007/BFb0064579},
}

@misc{dagnino_quantitative_2022,
	title = {Quantitative {Equality} in {Substructural} {Logic} via {Lipschitz} {Doctrines}},
	url = {http://arxiv.org/abs/2110.05388},
	doi = {10.48550/arXiv.2110.05388},
	urldate = {2023-03-01},
	publisher = {arXiv},
	author = {Dagnino, Francesco and Pasquali, Fabio},
	month = nov,
	year = {2022},
	note = {arXiv:2110.05388 [cs, math]},
}

@article{shulmanEnrichedIndexedCategories2013,
  title = {Enriched {{Indexed Categories}}},
  author = {Shulman, Michael},
  date = {2013},
  journaltitle = {Theory \& Applications of Categories},
  volume = {28},
  number = {21},
  url = {http://www.tac.mta.ca/tac/volumes/28/21/28-21abs.html},
  urldate = {2024-03-26}
}

@book{meyer-nieberg_banach_1991,
	series = {Universitext},
	title = {Banach lattices},
	isbn = {3-540-54201-9},
	url = {https://books.google.com/books?hl=en&lr=&id=2TPvCAAAQBAJ&oi=fnd&pg=PA1&dq=banach+lattices&ots=2VwdGWNxPi&sig=gRyjdQnom8SXtClFsohLHhvbsC4},
	urldate = {2024-03-26},
	publisher = {Springer-Verlag},
	author = {Meyer-Nieberg, Peter},
	year = {1991},
}

@inproceedings{glorot_deep_2011,
	title = {Deep sparse rectifier neural networks},
	url = {http://proceedings.mlr.press/v15/glorot11a},
	urldate = {2024-04-16},
	booktitle = {Proceedings of the fourteenth international conference on artificial intelligence and statistics},
	publisher = {JMLR Workshop and Conference Proceedings},
	author = {Glorot, Xavier and Bordes, Antoine and Bengio, Yoshua},
	year = {2011},
	pages = {315--323},
}

@article{nielsen_guaranteed_2016,
	title = {Guaranteed bounds on the {Kullback}-{Leibler} divergence of univariate mixtures using piecewise log-sum-exp inequalities},
	volume = {18},
	issn = {1099-4300},
	url = {http://arxiv.org/abs/1606.05850},
	doi = {10.3390/e18120442},
	language = {en},
	number = {12},
	urldate = {2024-04-16},
	journal = {Entropy},
	author = {Nielsen, Frank and Sun, Ke},
	month = dec,
	year = {2016},
	note = {arXiv:1606.05850 [cs, math, stat]},
	pages = {442},
}

@article{lawvere_adjointness_1969,
	title = {Adjointness in foundations},
	url = {https://www.jstor.org/stable/42969800?casa_token=GL45s_dCKqQAAAAA:SV1lSO9MTQgPgPiwpYRwbX2Dt8agvDzsifxWybePlKcVBm7MqgcDrfWueuE94LKuuKEutiHxK1tbA6NUup1AKFMCyruGhGpg00AD8_I4c4KNg760knXF},
	urldate = {2024-03-22},
	journal = {Dialectica},
	author = {Lawvere, F. William},
	year = {1969},
	note = {Publisher: JSTOR},
	pages = {281--296},
}

@article{lawvere_metric_1973,
	title = {Metric spaces, generalized logic, and closed categories},
	volume = {43},
	journal = {Rendiconti del seminario matématico e fisico di Milano},
	author = {Lawvere, F. William},
	year = {1973},
	note = {Publisher: Springer},
	pages = {135--166},
}

@article{grandis_categories_2007,
	title = {Categories, norms and weights},
	volume = {2},
	language = {en},
	journal = {Journal of Homotopy and Related Structures},
	author = {Grandis, Marco},
	year = {2007},
}

@inproceedings{dagnino_logical_2022,
	address = {New York, NY, USA},
	series = {{LICS} '22},
	title = {Logical {Foundations} of {Quantitative} {Equality}},
	isbn = {978-1-4503-9351-5},
	url = {https://doi.org/10.1145/3531130.3533337},
	doi = {10.1145/3531130.3533337},
	urldate = {2023-03-03},
	booktitle = {Proceedings of the 37th {Annual} {ACM}/{IEEE} {Symposium} on {Logic} in {Computer} {Science}},
	publisher = {Association for Computing Machinery},
	author = {Dagnino, Francesco and Pasquali, Fabio},
	month = aug,
	year = {2022},
	pages = {1--13},
}

@article{shulman_affine_2022,
	title = {Affine logic for constructive mathematics},
	volume = {28},
	issn = {1079-8986, 1943-5894},
	url = {https://www.cambridge.org/core/journals/bulletin-of-symbolic-logic/article/abs/affine-logic-for-constructive-mathematics/F710EB6032D94D936B5F144EBE124B58},
	doi = {10.1017/bsl.2022.28},
	language = {en},
	number = {3},
	urldate = {2023-07-25},
	journal = {Bulletin of Symbolic Logic},
	author = {Shulman, Michael},
	month = sep,
	year = {2022},
	note = {Publisher: Cambridge University Press},
	pages = {327--386},
}

@article{zadeh_fuzzy_1965,
	title = {Fuzzy sets},
	volume = {8},
	url = {https://www.sciencedirect.com/science/article/pii/S001999586590241X},
	number = {3},
	urldate = {2024-05-16},
	journal = {Information and control},
	author = {Zadeh, Lotfi Asker},
	year = {1965},
	note = {Publisher: Elsevier},
	pages = {338--353},
}

@article{zadeh_fuzzy_1988,
	title = {Fuzzy logic},
	volume = {21},
	url = {https://ieeexplore.ieee.org/abstract/document/53/},
	number = {4},
	urldate = {2024-05-16},
	journal = {Computer},
	author = {Zadeh, Lotfi Asker},
	year = {1988},
	note = {Publisher: IEEE},
	pages = {83--93},
}

@article{barr_variable_1985,
	title = {Variable set theory},
	url = {https://www.math.mcgill.ca/barr/papers/vst.pdf},
	urldate = {2024-05-10},
	journal = {Unpublished manuscript},
	author = {Barr, Michael and Mclarty, Colin and Wells, Charles},
	year = {1985},
}

@article{Figueroa2022,
  title = {A topos for continuous logic},
  author = {Figueroa, Daniel and van den Berg, Benno},
  journal = {Theory and Applications of Categories},
  volume = {38},
  number = {28},
  pages = {1108-1135},
  year = {2022},
  url = {http://www.tac.mta.ca/tac/volumes/38/28/38-28.pdf},
  MSC2020 = {03C66, 03G30, 18F10},
  note = {Published 2022-09-07}
}

@book{leinster_entropy_2021,
	title = {Entropy and diversity: the axiomatic approach},
	isbn = {978-1-108-96355-8},
	shorttitle = {Entropy and diversity},
	url = {https://books.google.com/books?hl=en&lr=&id=pI4qEAAAQBAJ&oi=fnd&pg=PP1&dq=Entropy+and+Diversity:+The+Axiomatic+Approach+leinster&ots=1Quru1FEVW&sig=RNC5qFY1AqkYkMd5lM_8PgQpEos},
	urldate = {2024-06-03},
	publisher = {Cambridge university press},
	author = {Leinster, Tom},
	year = {2021},
}

@book{boole_mathematical_1847,
	title = {The {Mathematical} {Analysis} of {Logic}},
	copyright = {Public domain in the USA.},
	url = {https://www.gutenberg.org/ebooks/36884},
	language = {English},
	urldate = {2024-06-03},
	publisher = {Barclay, Macmillan and Co.},
	author = {Boole, George},
	year = {1847},
}

@article{lukasiewicz_o_1920,
	title = {O logice trójwartościowej},
	volume = {5},
	journal = {Ruch Filozoficzny},
	author = {Łukasiewicz, Jan},
	year = {1920},
	pages = {170--171},
}

@article{scott1967proof,
  title={A Proof of the Independence of the Continuum Hypothesis},
  author={Scott, Dana},
  journal={Mathematical systems theory},
  volume={1},
  number={2},
  pages={89--111},
  year={1967},
  publisher={Springer-Verlag},
  url={https://www.karlin.mff.cuni.cz/~krajicek/scott67.pdf}
}

@article{cox_probability_1946,
	title = {Probability, {Frequency} and {Reasonable} {Expectation}},
	volume = {14},
	issn = {0002-9505},
	url = {https://doi.org/10.1119/1.1990764},
	doi = {10.1119/1.1990764},
	number = {1},
	urldate = {2024-06-03},
	journal = {American Journal of Physics},
	author = {Cox, R. T.},
	month = jan,
	year = {1946},
	pages = {1--13},
}

@article{fritz_synthetic_2020,
	title = {A synthetic approach to {Markov} kernels, conditional independence and theorems on sufficient statistics},
	volume = {370},
	issn = {00018708},
	url = {http://arxiv.org/abs/1908.07021},
	doi = {10.1016/j.aim.2020.107239},
	urldate = {2023-09-27},
	journal = {Advances in Mathematics},
	author = {Fritz, Tobias},
	month = aug,
	year = {2020},
	note = {arXiv:1908.07021 [cs, math, stat]},
	pages = {107239},
}

@article{jacobs_mathematics_2019,
	title = {The mathematics of changing one's mind, via {Jeffrey}'s or via {Pearl}'s update rule},
	volume = {65},
	issn = {1076-9757},
	url = {https://doi.org/10.1613/jair.1.11349},
	doi = {10.1613/jair.1.11349},
	number = {1},
	urldate = {2023-04-16},
	journal = {Journal of Artificial Intelligence Research},
	author = {Jacobs, Bart},
	month = may,
	year = {2019},
	pages = {783--806},
}

@misc{di_lavore_evidential_2023,
	title = {Evidential {Decision} {Theory} via {Partial} {Markov} {Categories}},
	url = {http://arxiv.org/abs/2301.12989},
	language = {en},
	urldate = {2024-03-13},
	publisher = {arXiv},
	author = {Di Lavore, Elena and Román, Mario},
	month = apr,
	year = {2023},
	note = {arXiv:2301.12989 [cs, math]},
}

@inproceedings{kozen_probabilistic_1983,
	address = {Not Known},
	title = {A probabilistic {PDL}},
	isbn = {978-0-89791-099-6},
	url = {http://portal.acm.org/citation.cfm?doid=800061.808758},
	doi = {10.1145/800061.808758},
	language = {en},
	urldate = {2024-06-07},
	booktitle = {Proceedings of the fifteenth annual {ACM} symposium on {Theory} of computing  - {STOC} '83},
	publisher = {ACM Press},
	author = {Kozen, Dexter},
	year = {1983},
	pages = {291--297},
}

@article{diff1,
  author       = {Samuel Teuber and
                  Stefan Mitsch and
                  Andr{\'{e}} Platzer},
  title        = {Provably Safe Neural Network Controllers via Differential Dynamic
                  Logic},
  journal      = {CoRR},
  volume       = {abs/2402.10998},
  year         = {2024},
  url          = {https://doi.org/10.48550/arXiv.2402.10998},
  doi          = {10.48550/ARXIV.2402.10998},
  eprinttype    = {arXiv},
  eprint       = {2402.10998},
  bibsource    = {dblp computer science bibliography, https://dblp.org}
}

@inproceedings{diff2,
  author       = {Natalia Slusarz and
                  Ekaterina Komendantskaya and
                  Matthew L. Daggitt and
                  Robert J. Stewart and
                  Kathrin Stark},
  editor       = {Ruzica Piskac and
                  Andrei Voronkov},
  title        = {Logic of Differentiable Logics: Towards a Uniform Semantics of {DL}},
  booktitle    = {{LPAR} 2023: Proceedings of 24th International Conference on Logic
                  for Programming, Artificial Intelligence and Reasoning, Manizales,
                  Colombia, 4-9th June 2023},
  series       = {EPiC Series in Computing},
  volume       = {94},
  pages        = {473--493},
  publisher    = {EasyChair},
  year         = {2023},
  url          = {https://doi.org/10.29007/c1nt},
  doi          = {10.29007/C1NT},
  bibsource    = {dblp computer science bibliography, https://dblp.org}
}

@misc{perrone_kan_2021,
	title = {Kan extensions are partial colimits},
	url = {http://arxiv.org/abs/2101.04531},
	language = {en},
	urldate = {2024-03-22},
	publisher = {arXiv},
	author = {Perrone, Paolo and Tholen, Walter},
	month = feb,
	year = {2021},
	note = {arXiv:2101.04531 [math]},
}

@article{bashirovMultiplicativeCalculusIts2008,
  title = {Multiplicative Calculus and Its Applications},
  author = {Bashirov, Agamirza E. and Kurp\i nar, Emine M\i s\i rl\i{} and \"Ozyap\i c\i, Ali},
  date = {2008-01-01},
  journaltitle = {Journal of Mathematical Analysis and Applications},
  shortjournal = {Journal of Mathematical Analysis and Applications},
  volume = {337},
  number = {1},
  pages = {36--48},
  issn = {0022-247X},
  doi = {10.1016/j.jmaa.2007.03.081},
  url = {https://www.sciencedirect.com/science/article/pii/S0022247X07003824},
  urldate = {2024-05-21}
}

\end{document}